\documentclass{amsart}

\RequirePackage{amsmath}
\RequirePackage{amsthm}
\RequirePackage{amsfonts}
\RequirePackage{amssymb}
\RequirePackage{enumerate} 
\RequirePackage{cite}
\RequirePackage{hyperref}

\newtheorem{theorem}{Theorem}[section]
\newtheorem{proposition}[theorem]{Proposition}
\newtheorem{corollary}[theorem]{Corollary}
\newtheorem{lemma}[theorem]{Lemma}
\theoremstyle{definition}
\newtheorem{definition}[theorem]{Definition}
\newtheorem{remark}[theorem]{Remark}

\newcommand{\C}{\mathbb{C}}

\newcommand{\R}{\mathbb{R}}
\newcommand{\D}{\mathbb{D}}
\newcommand{\N}{\mathbb{N}}

\renewcommand{\H}{\mathbb{H}}

\newcommand{\abs}[1]{\left| #1 \right|}

\title{Characterization of finite shift via Herglotz's representation}

\date{\today}
\thanks{This research was supported by Ministerio de Innovaci\'on y Ciencia, Spain, project PID2022-136320NB-I00 and by Ministerio de Universidades, Spain, through the action Ayuda del Programa de Formaci\'on de Profesorado Universitario, reference FPU21/00258.}
\author[F. J. Cruz-Zamorano]{Francisco J. Cruz-Zamorano}
\address{Departamento de Matem\'atica Aplicada II and IMUS, Escuela T\'ecnica Superior de Ingenier\'ia, Universidad de Sevilla,
	Camino de los Descubrimientos, s/n 41092, Sevilla, Spain}
\email{fcruz4@us.es}

\begin{document}
\begin{abstract}
A complete characterization of parabolic self-maps of finite shift is given in terms of their Herglotz's representation. This improves a previous result due to Contreras, D\'iaz-Madrigal, and Pommerenke. We also derive some consequences for the rate of convergence of these functions to their Denjoy-Wolff point, improving a related result of Kourou, Theodosiadis, and Zarvalis for the continuous setting.
\end{abstract}
\maketitle
\section{Introduction}
In the field of Discrete Complex Dynamics of the upper half-plane $\H = \{z = x+iy \in \C : y > 0\}$, the main object of study is the behaviour of the sequence of iterates of holomorphic self-maps $f \colon \H \to \H$, that is, the sequence of functions given by $f^{n+1} = f \circ f^n$, $n \in \N$, where we define $f^0 = \mathrm{Id}_{\H}$.

In this paper we focus on non-elliptic self-maps $f$, that is, those with no fixed points on $\H$. The dynamical behaviour of non-elliptic functions is described in the following seminal result:
\begin{theorem}[Denjoy-Wolff]
\cite[Theorem 3.2.1]{AbateBook}
\label{thm:DW}
Let $f \colon \H \to \H$ be a non-elliptic holomorphic map. Then, there exists $\tau \in \R \cup \{\infty\}$ such that $f^n \to \tau$ uniformly on compact sets of $\H$.
\end{theorem}
The point $\tau$ is commonly known as the Denjoy-Wolff point of $f$, and it acts as a global attractor of its dynamic.

Up to a conjugation with an isomorphism of $\H$, one can always suppose that the Denjoy-Wolff point of a non-elliptic self-map is at infinity. Indeed, this might be useful to find different non-elliptic functions whose behaviour varies. To this extent, we note that under this assumption one has \cite[Corollary 2.5.5]{AbateBook}
$$\angle\lim_{z \to \infty}f(z) = \infty, \quad \angle\lim_{z \to \infty}\dfrac{f(z)}{z} = \lambda \in [1,+\infty).$$
If $\lambda > 1$, $f$ is known as hyperbolic. This work deals with parabolic functions, that is, those for which $\lambda = 1$.

The study of parabolic functions will be done through a representation result which is due to Herglotz \cite[Theorem 6.2.1]{Aaronson}: every holomorphic function $f \colon \H \to \H$ can uniquely be written as
$$f(z) = \alpha z + \beta + \int_{\R}\dfrac{1+tz}{t-z}d\mu(t), \quad z \in \H,$$
where $\alpha \geq 0$, $\beta \in \R$ and $\mu$ is a positive finite measure on $\R$.

Using \cite[Chapter 5, Lemma 2]{DM}, one can deduce the following:
\begin{theorem}
\label{thm:Herglotz}
A holomorphic self-map $f \colon \H \to \H$ is parabolic with Denjoy-Wolff point at infinity if and only if it can be written as
\begin{equation}
\label{eq:Herglotz}
f(z) = z + \beta + \int_{\R}\dfrac{1+tz}{t-z}d\mu(t), \quad z \in \H,
\end{equation}
where $\beta \in \R$ and $\mu$ is a positive finite measure on $\R$, both of them not simultaneously null.
\end{theorem}

Following the latter result, every property of $f$ should be mirrored on the corresponding representing parameters $\beta$ and $\mu$: this is the idea behind this paper. In particular, we are interested in characterizing the following property:
\begin{definition}
\label{def:shift}
Let $f \colon \H \to \H$ be a holomorphic self-map which is parabolic and whose Denjoy-Wolff point is at infinity. The self-map $f$ is said of finite (respectively, infinite) shift if there exists an initial point $z_0 \in \H$ for which the orbit $z_n = x_n+iy_n = f^n(z_0)$ is such that $\sup_n y_n < +\infty$ (respectively, $\sup_n y_n = +\infty$).
\end{definition}
\begin{remark}
\cite[Proposition 3.2]{CDP}
The latter definition does not depend on the initial point $z_0\in \H$.
\end{remark}

Functions of finite shift have attracted much attention in the area. For example, a characterization of this property in terms of angular regularity of the functions is given in \cite[Theorem 4.1]{CDP}. Poggi-Corradini also used a related property to work with backward orbits; see \cite{PC_BI}. There are some references in the continuous version of Complex Dynamics (that is, the theory of holomorphic semigroups) as well; see \cite{BCDM} for an introduction to this topic. For example, in \cite[Theorem 3]{DB}, \cite{NK}, and \cite[Theorem 1.1]{DC}, the authors characterize parabolic semigroups of finite shift in terms of different geometric quantities.

Beyond functions of finite shift, in \cite{HS}, Theorem \ref{thm:Herglotz} has also proved to be useful to characterize the hyperbolic step of parabolic functions (see Section \ref{sec:preliminaries} for a definition). Indeed, many of the tools of this work have been previously derived there.

The main result of this work is the following, which completely characterizes self-maps of finite shift:
\begin{theorem}
\label{thm:shift}
Let $f \colon \H \to \H$ be a parabolic map with its Denjoy-Wolff point at infinity. If $f$ is of finite shift, then $\int_{\R}\abs{t}d\mu(t) < + \infty$.

Moreover, $f$ is of finite shift if and only if one of the following holds:
\begin{enumerate}[\hspace{0.5cm}\normalfont(i)]
\item $\int_{(-\infty,0)}\abs{t}d\mu(t) < +\infty$, $\int_{(0,+\infty)}t^2d\mu(t) < + \infty$, and $\beta > \int_{\R}td\mu(t)$.
\item $\int_{(-\infty,0)}t^2d\mu(t) < +\infty$, $\int_{(0,+\infty)}\abs{t}d\mu(t) < + \infty$, and $\beta < \int_{\R}td\mu(t)$.
\end{enumerate}
\end{theorem}

The proof of the first part of this result is given in Subsection \ref{subsec:proof1}. Then, in Subsection \ref{subsec:characterization}, the second part is deduced as a consequence of different results in \cite{CDP,HS}. Notice that this result improves the description of functions of finite shift given in \cite[Theorem 4.1]{CDP}, as noted in Subsection \ref{subsec:remarks}.

As a consequence of Theorem \ref{thm:shift}, we obtain the following result:
\begin{corollary}
\label{cor:convergence}
Let $f \colon \H \to \H$ be a parabolic map of finite shift with its Denjoy-Wolff point at infinity. Then
$$\dfrac{f^n}{n} \to \beta - \int_{\R}td\mu(t), \quad n \to +\infty,$$
uniformly on compact sets of $\H$.
\end{corollary}

We end this paper with a discussion on the rate of convergence of parabolic functions of finite shift to its Denjoy-Wolff point. This issue has been extensively discussed in \cite{DB-Rate-H,DB-Rate-P,BCDM-Rate} under the general setting of continuous dynamics. In the specific setting of parabolic semigroups of finite shift, better estimates for the rate of convergence have been obtained in \cite{KTZ}.

To introduce this topic in detail, we should move from the upper half-plane $\H$ into the unit disk $\D$. Using a M\"obius map $S \colon \D \to \H$, this can be done by relating any non-elliptic self-map $g \colon \D \to \D$ with a self-map $f = S \circ g \circ S^{-1} \colon \H \to \H$ whose Denjoy-Wolff point is infinity.

Definition \ref{def:shift} can be translated into the setting a self-map $g$ of the unit disk, meaning that $g$ is of finite shift if there exists some orbit of $g$ and a horocycle (see \cite[Definition 2.1.1]{AbateBook}) of center $\tau$ in which the orbit does not enter. If this is not the case, $g$ is said of infinite shift. Moreover, Theorem \ref{thm:DW} also remains valid. The point $\tau = S^{-1}(\infty) \in \partial\D$ is called the Denjoy-Wolff of $g$ and it follows that $g^n \to \tau$ as $n \to +\infty$ uniformly on compact sets of $\D$. By discussing the rate of convergence to the Denjoy-Wolff point we mean to estimate how fast the latter convergence is. To this extent, we prove the next result in Section \ref{sec:rate}:

\begin{theorem}
\label{thm:rate}
Let $g \colon \D \to \D$ be a parabolic map of finite shift whose Denjoy-Wolff point is $\tau \in \partial\D$. Then, there exists $C = C(g) > 0$ such that
$$n\abs{g^n-\tau} \to C, \quad n \to +\infty,$$
uniformly on compact sets of $\D$.
\end{theorem}

As we will see at the end, this result can be used to improve the rate of convergence given in \cite[Theorem 1.1]{KTZ}.

\textbf{Acknowledgement.} I would like to thank E. K. Theodosiadis and K. Zarvalis for their great hospitality and our fruitful conversations during my stay at the Aristotle University of Thessaloniki.

\section{Preliminaries: parabolic dynamics.}
\label{sec:preliminaries}
Parabolic self-maps are further classified in terms of their hyperbolic step. To define it, let $\rho$ be the pseudo-hyperbolic distance on $\H$. For any initial point $z_0 \in \H$, the Scharwz-Pick Lemma assures that the sequence $\rho(f^{n+1}(z_0),f^n(z_0))$ is non-increasing with respect to $n \in \N$. Therefore, it must converge. Indeed, if it converges to zero, then it does so for every initial point \cite[Corollary 4.6.9.(i)]{AbateBook}. Thus, the following definition is given:
\begin{definition}
A parabolic self-map $f \colon \H \to \H$ is said of zero hyperbolic step if $\rho(f^{n+1}(z_0),f^n(z_0)) \to 0$ as $n \to + \infty$ for some (every) $z_0 \in \H$. In another case, it is said of positive hyperbolic step.
\end{definition}

Determining the hyperbolic step of a given function can be a subtle task. However, some theoretical results are known:
\begin{theorem}
\cite{PommerenkeHalfPlane}
\label{thm:Pommerenke}
Let $z_0 \in \H$ be an initial point, and consider its orbit $z_n = x_n+iy_n = f^n(z_0)$, $n \in \N$. The following limit
$$b = \lim_{n \to \infty}\dfrac{x_{n+1}-x_n}{y_n} \in \R$$
exists. Furthermore, $b = 0$ if and only if $f$ is of zero hyperbolic step. Also, $z_{n+1}/z_n \to 1$ (see \cite[Eq. (3.16)]{PommerenkeHalfPlane}), and $y_{n+1}/y_n \to 1$ (see \cite[Eq. (3.17)]{PommerenkeHalfPlane}).
\end{theorem}
\begin{theorem}
\label{thm:shift-hs}
\cite[Proposition 3.3]{CDP}
Parabolic self-maps of finite shift are also of positive hyperbolic step.
\end{theorem}
\begin{remark}
There are also self-maps of infinite shift and positive hyperbolic step. For examples of this with the aid of the representation of Theorem \ref{thm:Herglotz}, we refer to \cite[Theorem 1.6]{HS}.
\end{remark}

\section{Proof of Theorem \ref{thm:shift} and Corollary \ref{cor:convergence}}
\label{sec:proof}
In this section, from now on, the function $f \colon \H \to \H$ will always be a holomorphic self-map which is parabolic and whose Denjoy-Wolff point is at infinity. Indeed, with the aid of Theorem \ref{thm:Herglotz}, we will identify $f$ with $(\beta,\mu)$.

\subsection{Proof of the first part of Theorem \ref{thm:shift}.}
\label{subsec:proof1}
Assume that $f$ is of finite shift and let us show that $\int_{\R}\abs{t}d\mu(t) < +\infty$. Notice that, by Theorem \ref{thm:shift-hs}, $f$ is of positive hyperbolic step. In particular, Theorem \ref{thm:Pommerenke} assures that for a given initial point $z_0 \in \H$, the orbit $z_n = x_n+iy_n = f^n(z_0)$ is such that
\begin{equation}
\label{eq:Pommerenke}
\lim_{n \to \infty}\dfrac{x_{n+1}-x_n}{y_n} = L \in \R \setminus \{0\}.
\end{equation}
But, as $f$ is of finite shift, there must exists $Y > 0$ such that 
\begin{equation}
\label{eq:Y}
\lim_{n \to + \infty}y_n = Y.
\end{equation}
Thus, $x_{n+1}-x_n \to LY =: \Delta \neq 0$. In the following, we will suppose that $\Delta > 0$ (otherwise, the proof will apply with slight modifications). In that case, Stolz's lemma assures that
\begin{equation}
\label{eq:Delta}
\lim_{n \to + \infty} \dfrac{x_n}{n\Delta} = 1.
\end{equation}

Under this setting, we can show the following:
\begin{lemma}
\label{lemma:t2}
$\int_{(0,+\infty)}t^2d\mu(t) < +\infty$.
\begin{proof}
First of all, notice that
$$y_{n+1} = y_n\left(1+\dfrac{y_{n+1}-y_n}{y_n}\right) = y_0\prod_{k=0}^n\left(1+\dfrac{y_{k+1}-y_k}{y_k}\right), \quad n \in \N.$$
Thus, the following are equivalent:
\begin{enumerate}[\hspace{0.5cm}\normalfont(i)]
\item $y_n$ is bounded,
\item $\prod_{n=0}^{\infty}\left(1+\frac{y_{n+1}-y_n}{y_n}\right) < +\infty,$
\item $\sum_{n = 0}^{\infty}\frac{y_{n+1}-y_n}{y_n} < +\infty.$
\end{enumerate}
But, from \eqref{eq:Herglotz},
$$\sum_{n = 0}^{\infty}\dfrac{y_{n+1}-y_n}{y_n} = \sum_{n = 0}^{+\infty} \dfrac{\mathrm{Im}(f(z_n)-z_n)}{y_n}= \int_{\R}\left(\sum_{n = 0}^{\infty}\dfrac{1+t^2}{(t-x_n)^2+y_n^2}\right)d\mu(t).$$

Using \eqref{eq:Delta}, let $N \in \N$ be such that
$$\dfrac{\Delta}{2} \leq x_{n+1}-x_n \leq \dfrac{3\Delta}{2}, \quad n \geq N.$$
Thus, for every $t > x_N$ one can find $k = k(t) \geq N$ such that $\abs{t-x_k} \leq \Delta$. In that case, using also \eqref{eq:Y},
$$\sum_{n \in \N}\dfrac{1+t^2}{(t-x_n)^2+y_n^2} \geq \dfrac{1+t^2}{\Delta^2+y_k^2} \geq \dfrac{1+t^2}{\Delta^2+Y^2}, \quad t > x_N.$$
It follows that,
$$+\infty > \int_{\R}\left(\sum_{n = 0}^{\infty}\dfrac{1+t^2}{(t-x_n)^2+y_n^2}\right)d\mu(t) \geq \int_{(x_N,+\infty)}\dfrac{1+t^2}{\Delta^2+Y^2}d\mu(t).$$
Since $\mu$ is a positive finite measure, the result follows.
\end{proof}
\end{lemma}

To continue with the proof of the first part of Theorem \ref{thm:shift}, let us proceed by supposing that
\begin{equation}
\label{eq:l1}
\int_{\R}\abs{t}d\mu(t) = + \infty.
\end{equation}
In that case, by Lemma \ref{lemma:t2}, it has to be that
$$\int_{(-\infty,0)}\abs{t}d\mu(t) = +\infty.$$
Under this hypothesis, we will prove that $x_{n+1}-x_n \to +\infty$, which contradicts \eqref{eq:Delta}, therefore showing that \eqref{eq:l1} cannot hold. Thus, the result will follow.

To do so, write $f$ as
$$f(z) = z + \tilde{\beta} + p_1(z) + p_0(z), \quad z \in \H,$$
where
$$\omega(A) = \mu(A \cap (-\infty,0]), \quad \dfrac{d\nu}{d\mu}(t) = (1+t^2)\chi_{(0,+\infty)}(t),$$
and
$$\tilde{\beta} = \beta -\int_{(0,+\infty)}td\mu(t) \in \R, \quad p_1(z) = \int_{\R}\dfrac{1+tz}{t-z}d\omega(t), \quad p_0(z) = \int_{\R}\dfrac{d\nu(t)}{t-z}, \quad z \in \H.$$
Note that, by Lemma \ref{lemma:t2}, $\tilde{\beta}$ is well-defined and the measure $\nu$ is finite. Indeed, if $y_0 \geq 1$, notice that
$$\abs{\text{Re}(p_0(z_n))} \leq \abs{p_0(z_n)} \leq \int_{\R}\dfrac{d\nu(t)}{\abs{t-z_n}} \leq \dfrac{\nu(\R)}{y_n} \leq \nu(\R), \quad n \in \N.$$ 
Then, it suffices to prove that $\mathrm{Re}(p_1(z_n)) \to + \infty$. To see this, notice that $p_1$ is a well-defined function on $\C \setminus (-\infty,0]$. Indeed, write
$$\text{Re}(p_1(x+iy)) = \int_{\R}\dfrac{t-x+t^2x-t(x^2+y^2)}{(t-x)^2+y^2}d\omega(t),$$
and notice that
$$\dfrac{\partial}{\partial y}\left(\dfrac{t-x+t^2x-t(x^2+y^2)}{(t-x)^2+y^2}\right) = -\dfrac{2y(1+t^2)(t-x)}{((t-x)^2+y^2)^2} \geq 0,$$
whenever $t \leq 0 < x$, $y \geq 0$. Thus, by \eqref{eq:Delta}, find $N \in \N$ such that $x_n > 1$ for all $n \geq N$ and notice that
$$\mathrm{Re}(p_1(z_n)) \geq \mathrm{Re}(p_1(x_n)) = \int_{(-\infty,0]}\dfrac{1+tx_n}{t-x_n}d\omega(t),$$
which is well-defined. Since $x_n \to +\infty$, the proof follows from the fact that
\begin{equation}
\label{eq:infty}
\lim_{x \to +\infty}\int_{(-\infty,0]}\dfrac{1+tx}{t-x}d\omega(t) = +\infty.
\end{equation}
To prove \eqref{eq:infty}, consider $x > 0$ and notice that the integrand is positive if and only if $t < -1/x$. By Fatou's lemma, one has
$$\liminf_{x \to +\infty}\int_{(-\infty,-1/x)}\dfrac{1+tx}{t-x}d\omega(t) \geq \int_{(-\infty,0)}(-t)d\omega(t) = +\infty,$$
where we have used \eqref{eq:l1}. On the other hand,
$$\int_{[-1/x,0]}\abs{\dfrac{1+tx}{t-x}}d\omega(t) \leq \dfrac{1}{x}\omega([-1/x,0]) \to 0, \quad x \to + \infty.$$
Thus, \eqref{eq:infty} holds.

\subsection{Proof of the second part of Theorem \ref{thm:shift}.}
\label{subsec:characterization}
Recall the following result:
\begin{proposition}
\cite[Proposition 3.4]{HS}
Let $f \colon \H \to \H$ be a parabolic self-map with its Denjoy-Wolff point at infinity such that $\int_{\R}\abs{t}d\mu(t) < +\infty$. Then, $f$ is of positive hyperbolic step if and only if it is of finite shift.
\end{proposition}
Thus, the second part of Theorem \ref{thm:shift} follows from Subsection \ref{subsec:proof1} and \cite[Theorem 1.4]{HS}.

\subsection{Further remarks.}
\label{subsec:remarks}
Let us recall that functions of finite shift are also in connection with self-maps enjoying some regularity in an angular sense. To introduce these ideas, the following definition is given in \cite{CDP}:
\begin{definition}
Consider a holomorphic map $\phi \colon \D \to \C$ with a boundary fixed point $\tau \in \partial\D$, that is,
$$\angle\lim_{z \to \tau}\phi(z) = \tau.$$
Then, $\phi$ is said to be of angular-class of order $p \in \N$ at $\tau$, denoted as $\phi \in C_A^p(\tau)$, if there exist $c_1,\ldots,c_p$ and a holomorphic function $\gamma \colon \D \to \C$ such that
$$\phi(z) = \tau + \sum_{k = 1}^p\dfrac{c_k}{k!}(z-\tau)^k+\gamma(z), \quad z \in \D, \quad \angle\lim_{z \to \tau}\dfrac{\gamma(z)}{(z-\tau)^p} = 0.$$
\end{definition}
In particular, as a consequence of Julia-Wolff-Carath\'eodory Theorem \cite[Corollary 2.5.5]{AbateBook}, every non-elliptic function is of angular-class of first order at its Denjoy-Wolff point.

These angular-classes are intimately linked to functions of finite shift, as the following result suggests:
\begin{theorem}
\cite[Theorem 4.1]{CDP}
Let $\phi \colon \D \to \D$ be a parabolic self-map of positive hyperbolic step with Denjoy-Wolff point $\tau \in \partial\D$. Then, $\phi$ is of finite shift if and only if $\phi \in C_A^2(\tau)$.
\end{theorem}

Theorem \ref{thm:shift} can be seen as a refinement of the latter result. To see this, one should link \cite[Propositions 2.1]{CDP} and \cite[Lemma 3.3]{HS} to state the following:
\begin{proposition}
Let $f \colon \H \to \H$ be a parabolic self-map with its Denjoy-Wolff point at infinity such that $\int_{\R}\abs{t}d\mu(t) < +\infty$. Consider the self-map $\phi \colon \D \to \D$ which is conjugated to $f$ and has Denjoy-Wolff point $\tau \in \partial\D$. Then, $\phi \in C_A^2(\tau)$. 
\end{proposition}
However, the converse is not true. These remarks have already been noticed in \cite[p. 15]{HS}.

\subsection{Proof of Corollary \ref{cor:convergence}.}
We will prove that there exists a non-empty open set $U \subset \H$ such that $$\dfrac{f^n(z)}{n} \to \beta - \int_{\R}td\mu(t), \quad n \to +\infty,$$
for all $z \in U$. Then, the result will follow from the identity principle and the fact that the sequence $\{f^n/n\}$ is a normal family.

To do so, as $f$ is of finite shift, let us suppose that $f$ is as in Theorem \ref{thm:shift}.(i) (if it is as in Theorem \ref{thm:shift}.(ii), the proof will also work with slight modifications). Write $f$ as
$$f(z) = z + \beta - \int_{\R}td\mu(t)+p(z), \quad z \in \H.$$
By \cite[Lemma 3.9]{HS}, there exists $a,b > 0$ such that $p(f^n(z)) \to 0$ as $n \to +\infty$ if
$$z \in \Omega := \{z =x+iy \in \C : x > a, \, y > b\}.$$
In particular, if $z \in \Omega$, it follows that
$$\mathrm{Re}(f^{n+1}(z)) - \mathrm{Re}(f^n(z)) \to \beta - \int_{\R}d\mu(t), \quad n \to +\infty.$$
Applying Stolz's lemma, it follows that
$$\dfrac{\mathrm{Re}(f^n(z))}{n} \to \beta-\int_{\R}td\mu(t), \quad n \to +\infty.$$
Moreover, since $f$ is of finite shift, $\mathrm{Im}(f^n(z))$ is bounded, and so $\mathrm{Im}(f^n)/n \to 0$ as $n \to +\infty$. Then, the result follows.

\section{Application: rate of convergence on the unit disk}
\label{sec:rate}
This section aims to develop a proof of Theorem \ref{thm:rate} and to discuss its implications on the theory of continuous dynamics. To do so, recall that any self-map $g \colon \D \to \D$ whose Denjoy-Wolff point is $\tau \in \partial\D$ can be conjugated to a self-map $f = S \circ g \circ S^{-1} \colon \H \to \H$ whose Denjoy-Wolff point is infinity, where $S \colon \D \to \H$ is the M\"obius map given by
$$S(z) = i\dfrac{\tau+z}{\tau-z}, \quad z \in \D, \quad S^{-1}(w) = \tau\dfrac{w-i}{w+i}, \quad w \in \H.$$
From this, a computation shows:
\begin{lemma}
\label{lemma:gtof}
With the above notations,
$$g^n(z)-\tau= \dfrac{-2i\tau}{f^n(S(z))+i}, \quad z \in \D.$$
\end{lemma}

\subsection{Proof of Theorem \ref{thm:rate}.}
It follows from Lemma \ref{lemma:gtof} and Corollary \ref{cor:convergence}, having in mind that $f^n(S(z)) \to \infty$ as $n \to +\infty$ for all $z \in \D$.

\subsection{Continuous dynamics.}
The uniform convergence on compact sets given in Corollary \ref{cor:convergence} lets us relate the ``discrete'' rate of convergence in Theorem \ref{thm:rate} with its analogous ``continuous'' one. That is, the next result (which improves \cite[Theorem 1.1]{KTZ} and complements \cite[Theorem 1.(b)]{DB-Rate-P}) can also be derived:
\begin{theorem}
\label{thm:semigroup}
Let $\{\phi_t \colon \D \to \D\}$ be a parabolic semigroup of self-maps of finite shift whose Denjoy-Wolff point is $\tau \in \partial\D$. There exists $C = C(\phi_1) > 0$ such that
$$t\abs{\phi_t-\tau} \to C, \quad t \to +\infty,$$
uniformly on compact sets of $\D$.
\begin{proof}
Consider $g = \phi_1$. Notice that $g$ is a parabolic self-map of $\D$ which is of finite shift and has Denjoy-Wolff point $\tau \in \partial\D$. Moreover, for every $t > 0$, consider the decomposition $t = n+s$, where $n \in \N$ and $s \in [0,1)$. Notice that, given $z \in \D$, one has
$$\min_{z \in K}\abs{g^n(K) - \tau} \leq \abs{\phi_t(z) - \tau} \leq \max_{z \in K}\abs{g^n(K) - \tau},$$
where $K = \{\phi_s(z) : s \in [0,1]\} \subset \D$ is a compact set. Since $n = n(t)$ and $t/n \to 1$ as $t \to +\infty$, the result follows from Theorem \ref{thm:rate}.
\end{proof}
\end{theorem}
\begin{remark}
The constant in Theorem \ref{thm:rate} (respectively, in Theorem \ref{thm:semigroup}) is
$$C = \dfrac{2}{\abs{\beta - \int_{\R}td\mu(t)}},$$
where $\beta$ and $\mu$ are the representing parameters in terms of Theorem \ref{thm:Herglotz} of the map $f \colon \H \to \H$ which is conjugated to $g$ (respectively, to $\phi_1$) through the M\"obius map $S$ used in Lemma \ref{lemma:gtof}.
\end{remark}


\begin{thebibliography}{99}

\bibitem{Aaronson}
J. Aaronson.
\newblock {\em An Introduction to infinite ergodic theory}.
\newblock {Mathematical Surveys and Monographs} \textbf{50}. American Mathematical Society, 1997. 

\bibitem{AbateBook}
M. Abate.
\newblock {\em Holomorphic Dynamics on Hyperbolic Riemann Surfaces}, {\bf 89}.
\newblock Walter de Gruyter GmbH \& Co KG, (2022).

\bibitem{DB}
D. Betsakos.
\newblock {\em Geometric description of the classification of holomorphic semigroups}.
\newblock {Proc. Amer. Math. Soc.} {\bf 144} (2016), 1595--1604.

\bibitem{DB-Rate-H}
D. Betsakos.
\newblock {\em On the rate of convergence of hyperbolic semigroups of holomorphic functions}.
\newblock {Bull. Lond. Math. Soc.} {\bf 47} (2015), 493--500.

\bibitem{DB-Rate-P}
D. Betsakos.
\newblock {\em On the rate of convergence of parabolic semigroups of holomorphic functions}.
\newblock {Anal. Math. Phys.} {\bf 5} (2015), 207--216.

\bibitem{BCDM-Rate}
D. Betsakos, M. D. Contreras, and S. D\'iaz-Madrigal.
\newblock {\em On the rate of convergence of semigroups of holomorphic functions at the Denjoy-Wolff point}.
\newblock {Rev. Mat. Iberoam.} {\bf 36} (2020), 1659--1686.

\bibitem{BCDM} 
F. Bracci, M. D. Contreras, and S. D\'iaz-Madrigal. 
\newblock {\it Continuous semigroups of holomorphic self-maps of the unit disc}.
\newblock Springer Monographs in Mathematics. Springer, (2020).

\bibitem{HS}
M. D. Contreras, F. J. Cruz-Zamorano, and L. Rodr\'iguez-Piazza.
\newblock {\em Characterization of the hyperbolic step of parabolic functions}.
\newblock {Submitted.} Available in arXiv:2309.00402.

\bibitem{CDP}
M. D. Contreras, S. D\'iaz-Madrigal, and Ch. Pommerenke.
\newblock {\em Second angular derivatives and parabolic iteration in the unit disk}.
\newblock {Trans. Amer. Math. Soc.}, {\bf 362} (2010), 357--388.

\bibitem{DC}
D. Cordella.
\newblock {\em Holomorphic semigroups of finite shift in the unit disc}.
\newblock {J. Math. Anal. Appl.}, {\bf 513} (2022), 17pp.

\bibitem{DM} C. I. Doering and R. Mañ\'e.
\newblock {\em The dynamics of inner functions}. 
\newblock {Ensaios Matemáticos}, {\bf 3}. Sociedade Brasileira de Matemática, 1991.

\bibitem{NK}
N. Karamanlis.
\newblock {Angular derivatives and semigroups of holomorphic functions.}.
\newblock {Illinois J. Math.} {\bf 63} (2019), 403--424.

\bibitem{KTZ}
M. Kourou, E. K. Theodosiadis, and K. Zarvalis.
\newblock {\em Rate of convergence for holomorphic semigroups of finite shift}
\newblock {Submitted}. Available at arXiv:2403.06883, (2024).

\bibitem{PC_BI}
P. Poggi-Corradini.
\newblock {\em Backward-iteration sequences with bounded hyperbolic steps for analytic self-maps of the disk}.
\newblock {Rev. Mat. Iberoamericana}, {\bf 19} (2003), 943--970.

\bibitem{PommerenkeHalfPlane}
Ch. Pommerenke.
\newblock {\em On the iteration of analytic functions in a halfplane, {I}}.
\newblock {J. London Math. Soc. (2)}, {\bf 19} (1979), 439--447.
    
\end{thebibliography}
\end{document}